\input amstex
\documentstyle{amsppt}
\magnification1200
\hsize 6.0 true in
\vsize 8.5 true in
\parskip=\medskipamount
\NoBlackBoxes
\NoRunningHeads
\TagsOnRight

\loadeusm
\loadbold

\topmatter

\title An elementary approach to character sums over multiplicative subgroups \endtitle
\author Ke Gong \endauthor
\address Department of Mathematics, Henan University, Kaifeng, Henan 475004, P.R. China \endaddress
\email kg{\@}henu.edu.cn \endemail

\thanks
\indent This work was supported by the Natural Science Foundation of China (Grant No. 11126150), the Natural Science Foundation of the Education Department of Henan Province (Grant No. 2011A110003), and a Qu\'ebec Merit Scholarship (2011--2012).
\endthanks

\abstract Let $p$ be an odd prime. Using I.~M.~Vinogradov's bilinear estimate, we present an elementary approach to estimate nontrivially the character sum
$$
\sum_{x\in H}\chi(x+a),\qquad a\in\Bbb F_p^*,
$$
where $H<\Bbb F_p^*$ is a multiplicative subgroup in finite prime field $\Bbb F_p$. Some interesting mean-value estimates are also provided.
\endabstract

\endtopmatter

\def\Z{{\Bbb Z}}
\def\F{{\Bbb F}}

\def\eps{\varepsilon}
\def\({\left(}
\def\){\right)}

\def\spmod#1{\,(\text{\rm mod\,}#1)}

\document

\noindent
{\bf 1. Introduction}

Let $p$ be an odd prime, and $H<\F_p^*$ be a multiplicative subgroup of the finite prime field $\F_p$. After his opening work on extremely short exponential sums $\sum_{x\in H}\psi(ax)$ with $\psi$ being the additive character of $\F_p$, Jean Bourgain posed the following problem concerning multiplicative character sums over shifted subgroup, see [C, Problem 5].

\proclaim{Problem 1~{\rm (J.~Bourgain)}} Obtain nontrivial bound on $\sum_{x\in H}\chi(x+a)$ for $H<\F_p^*$, $|H|\sim\sqrt{p}$, and $a\in\F_p^*$.
\endproclaim

Using I. M. Vinogradov's bilinear estimate for character sums, we shall present an elementary approach to Bourgain's character sum.

\proclaim{Theorem 2} For any $H<\F_p^*$, we have
$$
\max_{a\in\F_p^*}\left|\sum_{x\in H}\chi(x+a)\right|<p^{1/2}.
$$
\endproclaim

Thus, for any $\eps>0$ and $H<\F_p^*$ with $|H|>p^{1/2+\eps}$, we have
$$
\max_{a\in\F_p^*}\left|\sum_{x\in H}\chi(x+a)\right|<p^{-\varepsilon}|H|.
$$

We also obtain two mean-value estimates which suggest that an estimate for extremely short character sums may exist.

\bigskip
\noindent
{\bf 2. Vinogradov Lemma}

We recall that, first in 1930s (see [V43]) and then in his monograph [V76, Russian, p.\,88; English, pp.\,360--361], I. M. Vinogradov obtained the following bilinear estimate (up to a $\sqrt{2}$--factor in the upper bound) for character sums, which played a fundamental role in his studies on character sums over shifted primes.

\proclaim{Lemma 3~{\rm (I.~M.~Vinogradov)}} Let $p$ be an odd prime, $\gcd(a,p)=1$, $\chi\ne\chi_0\pmod p$,
$$
S = \sum_{x=0}^{p-1}\sum_{y=0}^{p-1}\xi(x)\eta(y)\chi(xy+a),
$$
$$
S' = \sum_{x=0}^{p-1}\sum_{y=0}^{p-1}\xi(x)\eta(y)\chi(xy(xy+a)).
$$
For any complex-valued functions $\xi$ and $\eta$ with
$$
\sum_{x=0}^{p-1}|\xi(x)|^2\le X,\qquad \sum_{y=0}^{p-1}|\eta(y)|^2\le Y,
$$
we have
$$
|S|\le\sqrt{pXY},\qquad |S'|\le\sqrt{pXY}.
$$
\endproclaim

We present here, for the sake of completeness, the proof due to Vinogradov [V81, Chap.\,V, Exercise 8, {\bf c}], where the Legendre symbol case was treated. Indeed, Lemma 3 is a counterpart of Vinogradov's bilinear estimate for exponential sums, see [V81, Chap.\,VI, Exercise 8, $\alpha$)].

\demo{Proof} It suffices to prove the first statement, since the second one is immediately if we take $\xi'(x)=\xi(x)\chi(x)$ and $\eta'(y)=\eta(y)\chi(y)$ as the weights.

Since
$$
\aligned
|S|^2
& \le \(\sum_{x=0}^{p-1}|\xi(x)|\cdot\left|\sum_{y=0}^{p-1}\eta(y)\chi(xy+a)\right|\)^2 \\
& \le X\sum_{x=0}^{p-1}\sum_{y=0}^{p-1}\sum_{y_1=0}^{p-1}\eta(y)\overline{\eta(y_1)}
\chi(xy+a)\overline\chi(xy_1+a) \\
& = X\sum_{y=0}^{p-1}\sum_{y_1=0}^{p-1}\eta(y)\overline{\eta(y_1)}S_{y,y_1},
\endaligned
$$
where
$$
S_{y,y_1} = \sum_{x=0}^{p-1}\chi(xy+a)\overline\chi(xy_1+a).
$$
Thus
$$
S_{y,y_1}=
\cases
                     p, & \text{if}~y=y_1=0; \\
                     0, & \text{only one of}~y~\text{and}~y_1~\text{equals}~0; \\
                   p-1, & \text{if}~y=y_1>0; \\ %% x\equiv -ay^{-1}\pmod p, the term is 0
-\chi\(\frac{y}{y_1}\), & \text{otherwise}.
\endcases
$$
The last equality is because
$$
\multline
S_{y,y_1} = \sum_{z=0}^{p-1} \chi\(\frac{y}{y_1}z+a\(1-\frac{y}{y_1}\)\)\overline\chi(z) \\
= \chi\(\frac{y}{y_1}\) \sum_{z=1}^{p-1}\chi(z)\overline\chi(z)\chi\(1+a\frac{y_1}{y}\(1-\frac{y}{y_1}\)z^*\)
= -\chi\(\frac{y}{y_1}\)
\endmultline
$$
by taking $xy_1+a\equiv z\pmod p$.
Therefore
$$
\aligned
|S|^2
& \le X \(p|\eta(0)|^2 + (p-1)\sum_{y=1}^{p-1}\eta(y)\overline{\eta(y)}
- \underset{y\ne y_1}\to{\sum_{y=1}^{p-1}\sum_{y_1=1}^{p-1}}\eta(y)\overline{\eta(y_1)}\chi\(\frac{y}{y_1}\)\) \\
& = X \(p|\eta(0)|^2 + p\sum_{y=1}^{p-1}\eta(y)\overline{\eta(y)}
- \sum_{y=1}^{p-1}\eta(y)\overline{\eta(y)}
- \underset{y\ne y_1}\to{\sum_{y=1}^{p-1}\sum_{y_1=1}^{p-1}}\eta(y)\overline{\eta(y_1)}\chi\(\frac{y}{y_1}\)\) \\
& = X \(p\sum_{y=0}^{p-1}\eta(y)\overline{\eta(y)}
- \sum_{y=1}^{p-1}\sum_{y_1=1}^{p-1}\eta(y)\overline{\eta(y_1)}\chi\(\frac{y}{y_1}\)\) \\
& = X \(p\sum_{y=0}^{p-1}|\eta(y)|^2 - \left|\sum_{y=1}^{p-1}\eta(y)\chi(y)\right|^2\)
\le pXY,
\endaligned
$$
which completes the proof. \qed
\enddemo

\bigskip
\noindent
{\bf 2. Proof of Theorem 2}

Here and below, we denote $A(\cdot)$ the indicator function for a subset $A$ of $\F_p$.

Indeed, once taking $\eta$ to be the indicator function of multiplicative subgroup $H<\F_p^*$ in the proof of Lemma~3, we have
$$
|S|^2\le |H|\(p|H|-\Big|\sum_{x\in H<\F_p^*}\chi(x)\Big|^2\).
$$
But we only need a weak upper bound for our use below.

Recall that $a\in\F_p^*$. We first write
$$
\sum_{x\in H}\chi(x+a) = \frac{1}{|H|}\sum_{x,y\in H}\chi(xy+a) =  \frac{1}{|H|}\sum_x\sum_y H(x)H(y)\chi(xy+a),
$$
then apply Lemma~1 directly to obtain
$$
\left|\sum_{x\in H}\chi(x+a)\right|\le \frac{1}{|H|}\sqrt{p|H|\cdot|H|} = \sqrt{p}. \tag 1
$$

\bigskip
\noindent
{\bf 3. Mean-value estimate, I}

We find that the following identity could be obtained from generalizing the classical results due to Vinogradov [V81, Chap. VII, Exercise 1], Davenport and Erd\H{o}s [DE, Lemma 1].

\proclaim{Theorem 4} For any subset $D\subset\F_p^*$, we have the identity
$$
\sum_{a\in\F_p}\left|\sum_{x\in D}\chi(x+a)\right|^2 = p|D|-|D|^2. \tag 2
$$
\endproclaim

\demo{Proof} Indeed,
$$
\aligned
\sum_{a\in\F_p}\left|\sum_{x\in D}\chi(x+a)\right|^2
& = \sum_{a\in\F_p}\sum_{x,\,y\in D}\chi(x+a)\overline\chi(y+a) \\
& = \sum_{a\in\F_p}\sum_{x\in D}|\chi(x+a)|^2 + \sum\Sb x,\,y\in D\\ x\ne y\endSb\sum_{a\in\F_p} \chi(x+a)\overline\chi(y+a) \\
& = (p-1)|D| - |D|(|D|-1) \\
& = p|D| - |D|^2.
\endaligned
$$
Note that the second to last equality is due to the fact
$$
\sum_{a\in\F_p} \chi(x+a)\overline\chi(y+a) = -1,
$$
which is a consequence of the observation that the congruence
$$
x+a\equiv z(y+a)\pmod p
$$
establishes a one-to-one correspondence between all $a$ with $a\not\equiv -y$ and all $z$ with $z\not\equiv 1$. \qed
\enddemo

\comment
Thus, we have
$$
\sum_{a\in\F_p^*}\left|\sum_{x\in D}\chi(x+a)\right|^2 = p|D| - |D|^2 - \left|\sum_{x\in D}\chi(x)\right|^2,
$$
and
$$
\max_{a\in\F_p^*}\left|\sum_{x\in D}\chi(x+a)\right|^2
\le p|D| - |D|^2 - \left|\sum_{x\in D}\chi(x)\right|^2.
$$
When $D=H<\F_p^*$, we have
$$
\max_{a\in\F_p^*}\left|\sum_{x\in H}\chi(x+a)\right|
\le \sqrt{ p|H| - |H|^2 - \left|\sum_{x\in H}\chi(x)\right|^2}.
$$
\endcomment

We remark that, (2) could be compared with its counterpart for exponential sums (see Konyagin [K, Lemma 2]). That is, for any subset $D\subset\Z_q$ ($q$ is a positive integer), there holds
$$
\sum_{a\in\Z_q\setminus\{0\}}\left|\sum_{x\in D}e_q(ax)\right|^2 = |D|(q-|D|).
$$

\bigskip
\noindent
{\bf 4. Mean-value estimate, II}

In this section we present another different type mean-value estimate. For $a\in\F_p^*$, we have
$$
\frac{1}{p-1}\sum_{\chi\spmod p}\left|\sum_{n\in H}\chi(n+a)\right|\le\sqrt{|H|}.
$$

\demo{Proof} Indeed,
$$
\aligned
\(\sum_{\chi\spmod p}\left|\sum_{n\in H}\chi(n+a)\right|\)^2
& \le (p-1)\sum_{m\in H}\sum_{n\in H}\sum_{\chi\spmod p}\chi\(\frac{m+a}{n+a}\) \\
& \le (p-1)((p-1)|H| + 0) \\
& = (p-1)^2|H|,
\endaligned
$$
which completes the proof. \qed
\enddemo

If $a=0$, we recall that Shkredov [S, p.\,607] has obtained $\sum_{\chi}\left|\sum_{n\in H}\chi(n)\right|\le p$. Indeed, we can even obtain an identity. Here we present a proof due to A.~Granville.

If $H$ is the subgroup of order $(p-1)/k$, then
$$
H(n) = \frac1k\sum_{\psi:\ \psi^k=\chi_0}\psi(n),
$$
so that
$$
\sum_{n} H(n)\chi(n) = \frac1k \sum_{\psi:\ \psi^k=\chi_0}\sum_{n}(\psi\chi)(n)
$$
and this equals $\frac{p-1}{k}$ if $\chi=\overline\psi$ for some $\psi$ satisfying $\psi^k= \chi_0$, and equals $0$ otherwise. Hence we see that
$$
\sum_{\chi}\Big|\sum_{n\in H}\chi(n)\Big|
= \sum_\chi\Big|\sum_{n} H(n)\chi(n)\Big|
= \sum_{\psi:\ \psi^k=\chi_0}\frac{p-1}{k} = p-1.
$$

\bigskip
\noindent
{\bf 5. Final remarks}

Firstly, estimate (1) can be obtained by Weil's estimate through expressing the indicator function of $H$ as $H(n) = \frac{1}{k}\sum_{\psi: \psi^k=\chi_0} \psi(n)$. However, our method is completely elementary.

Secondly, using the estimate for $S'$ in Lemma 3, we have for $a\in\F_p^*$
$$
\Big|\sum_{x\in H}\chi(x(x+a))\Big| = \frac{1}{|H|}\Big|\sum_x\sum_y H(x)H(y)\chi(xy(xy+a))\Big|
\le \sqrt{p},
$$
which provides for the nonlinear character sums the same upper bound as the linear case in Theorem 2.

We remark that I. E. Shparlinski has posed the following problem, which has immediate implications to polynomial factorization over finite fields.

\proclaim{Problem~5 {\rm (I. E. Shparlinski)}} Estimate nontrivially
$$
\sum_{x\in H}\chi((x+a)(x+b)),\qquad ab(a-b)\in\F_p^*
$$
for $H<\F_p^*$ and $|H|\sim\sqrt{p}$.
\endproclaim

Finally, we would like to pose the following problem, which could also be dealt directly by Weil's estimates for character sums with rational functions argument. But no elementary approach is known.

\proclaim{Problem 6} Estimate nontrivially the sums
$$
\sum_{x\in H}e\(\frac{kx+\ell x^*}{p}\),\qquad \sum_{x\in H\setminus\{-a\}}e\(\frac{k(x+a)^*}{p}\)
$$
for $H<\F_p^*$ with $|H|\sim\sqrt{p}$ and $k,\ell,a\in\F_p^*$.
\endproclaim

\bigskip
\noindent
{\bf Acknowledgement}

I would like to thank Professors Andrew Granville, Chaohua Jia and Igor Shparlinski for their interests and comments on this work. I am grateful to Professor Sergei Konyagin for sending his Russian paper.

%%%%%%%%%%% References %%%%%%%%%%%%%

\enddocument